\documentclass[letterpaper, preprint, paper,11pt]{AAS}	

\usepackage{bm}
\usepackage{kotex}
\usepackage{amsmath,amssymb,amsfonts}
\usepackage{subfigure}
\usepackage[colorlinks=true, pdfstartview=FitV, linkcolor=black, citecolor= black, urlcolor= black]{hyperref}
\usepackage{overcite}
\usepackage{footnpag}			      	
\usepackage{longtable,tabularx}
\usepackage{multirow}
\usepackage{algorithm,algcompatible}
\algnewcommand\INPUT{\item[\textbf{Input:}]}%
\algnewcommand\OUTPUT{\item[\textbf{Output:}]}%

\DeclareMathOperator*{\argmin}{argmin}

\PaperNumber{24-166}

\begin{document}

\title{Optimal Strip Attitude Command of Earth Observation Satellite using Differential Dynamic Programming}

\author{Seungyeop Han\thanks{PhD Student, Daniel Guggenheim School of Aerospace Engineering, Georgia Institute of Technology, Atlanta, GA,30332},\
Byeong-Un Jo\thanks{Assistant Professor, Department of Aerospace Engineering, Sejong University, Seoul 05006, Republic of Korea},\
 and Koki Ho\thanks{Dutton-Ducoffe Professor, Daniel Guggenheim School of Aerospace Engineering, Georgia Institute of Technology, Atlanta, GA, 30332}
}

\maketitle{}

\begin{abstract}
This paper addresses the optimal scan profile problem for strip imaging in an Earth observation satellite (EOS) equipped with a time-delay integration (TDI) camera. Modern TDI cameras can control image integration frequency during imaging operation, adding an additional degree of freedom (DOF) to the imaging operation. On the other hand, modern agile EOS is capable of imaging non-parallel ground targets, which require a substantial amount of angular velocity and angular acceleration during operation. We leverage this DOF to minimize various factors impacting image quality, such as angular velocity. Initially, we derive analytic expressions for angular velocity based on kinematic equations. These expressions are then used to formulate a constrained optimal control problem (OCP), which we solve using differential dynamic programming (DDP). We validate our approach through testing and comparison with reference methods across various practical scenarios. Simulation results demonstrate that our proposed method efficiently achieves near-optimal solutions without encountering non-convergence issues.
\end{abstract}
\section*{Notation}
\begin{tabular}{r l}
	$\textbf{\textit{a}}$  & geometric vector \\
    $\hat{\textbf{\textit{a}}}$ & unit vector of $\textbf{\textit{a}}$ \\
    $\textit{\textit{a}}$ & magnitude of $\textbf{\textit{a}}$ or scalar value \\
    $\textbf{\textit{a}}^A$ & $\textbf{\textit{a}}$ expressed in frame $A$ \\
    $\dot{a}$ & time derivative of scalar $a$ \\
    $^{A} \dot{\textbf{\textit{a}}}$ & time derivative of $\textbf{\textit{a}}$ with respect to frame $A$ \\
    $^{A} \dot{\textbf{\textit{a}}}^A$ & $^{A} \dot{\textbf{\textit{a}}}$ expressed in frame $A$ \\
\end{tabular} \\

\section{Introduction} 
High-resolution imagery and broader ground coverage have always been main areas of interest for Earth Observation Satellites (EOS). To fulfill these demands, many high-resolution EOS use time-delay integration (TDI) cameras, which enhance the signal-to-noise ratio and image quality by capturing multiple exposures of an object continuously and integrating each exposure. To achieve a high-quality image, it is necessary to synchronize the relative motion between the camera and the object. Recently, adjustable high scan rate TDI cameras have become popular for space applications. By adjusting the scan rate in real-time, the satellite system can capture a more general kind of ground scan with minimal image quality degradation.

Generally, there are two types of EOS: one consists of a large number of small-sized, low-resolution satellites forming a constellation to achieve global coverage, while the other comprises a small number of large-sized, high-resolution satellites designed to gather precise local data. The latter type of satellite is not easily mass-produced due to its production cost, size, and development period. Consequently, they are developed with agility to capture more images. For such satellites, the optimized imaging operation is necessary to maximize their utility within limited operational time.

Conventionally, and still in many cases today, the ground scan path of imaging operations is parallel to the EOS's ground trajectory. This type of imaging operation is known as 'strip' imaging because the captured images form a strip-like coverage. Two main reasons for maintaining a parallel ground scan path are the ease of motion synchronization and low angular velocity maneuvers. Some EOSs may lack the agility required to image non-parallel ground targets, or the resulting image quality may be poor, rendering such imaging operations impractical. However, with the availability of adjustable TDI cameras and agile EOS, strip imaging of non-parallel ground targets has been extensively studied, and this paper also focuses on this topic.

Related to non-parallel imaging operations, the work explained the time integral model of TDI cameras for a pre-determined attitude profile, and the research introduced an iterative method to find subsequent ground target points of non-parallel ground targets \cite{zhu2016integral, topal2016spaceborne}. The study described a method to compute an attitude specifically for north-directional imaging, and the paper introduced the imaging attitude for non-parallel ground targets, both with simplified geometric assumptions for both ground targets and satellite orbits \cite{jeon2018north, ye2020fuzzy}. The guidance algorithm with multiple piecewise ground target segments and the method to compute imaging attitudes for general targets and satellite orbits were also proposed, but both lack explanations of angular rate and acceleration \cite{ekinci2019guidance, du2022attitude}. The work was the first to design the scan rate of TDI cameras with consideration of both angular rate and acceleration, but it lacks an analytic expression for rate and acceleration and used simple parametric optimization rather than trajectory optimization to design the scan rate of TDI cameras \cite{qiu2020attitude}.

On the other hand, differential dynamic programming (DDP) is a classical but powerful optimization method for solving optimal control problems (OCP), utilizing Bellman's optimality principles and quadratic approximation to find local optimal solutions \cite{jacobson1970differential}. In recent years, it has gained popularity due to its fast convergence properties and advancements in research addressing general nonlinear constraints. The work has proposed DDP using a nominal trajectory with consideration of terminal constraints, and other research has handled box constraints on control inputs \cite{sun2014continuous,tassa2014control}. The active set method and the penalty method, combined with DDP, have been proposed to manage general nonlinear constraints \cite{xie2017differential,chen2019autonomous}. Additionally, the augmented Lagrangian method has been adopted for improved numerical stability and convergence speed \cite{plancher2017constrained,howell2019altro}.

This paper builds on the works, which seek analytical attitude profiles for parallel and general strip operations, respectively\cite{yoon2009analytical, han2017attitude}. Recently, the work explains the analytical expression of attitude, angular rate, and acceleration for general staring operations\cite{han2022analytical}. Utilizing these results, this paper derives analytical expressions of attitude profiles for general strip imaging operations. Additionally, the scan rate profile of the TDI camera is fully optimized using the DDP algorithm with the analytical attitude profile. Lastly, this paper demonstrates the proposed algorithm through numerical simulations under various strip scenarios.

\section{Analytic Expression of General Strip Attitude Profile}
\subsection{General Staring Attitude}
The objective of imaging attitude control is to align the sensing axis of a spacecraft (the boresight axis) with a target. The desired direction of the payload's sensing axis is determined by the relative motion between the target and spacecraft, which varies with time. Referring to Figure~\ref{fig_1}, the line of sight vector (LOS) $\boldsymbol{\rho}$, the LOS direction vector $\hat{\boldsymbol{\rho}}$, and the relative distance $\rho$ from the spacecraft to the target are defined as:

\begin{figure}[h]
\centering
\includegraphics[width=0.45\textwidth]{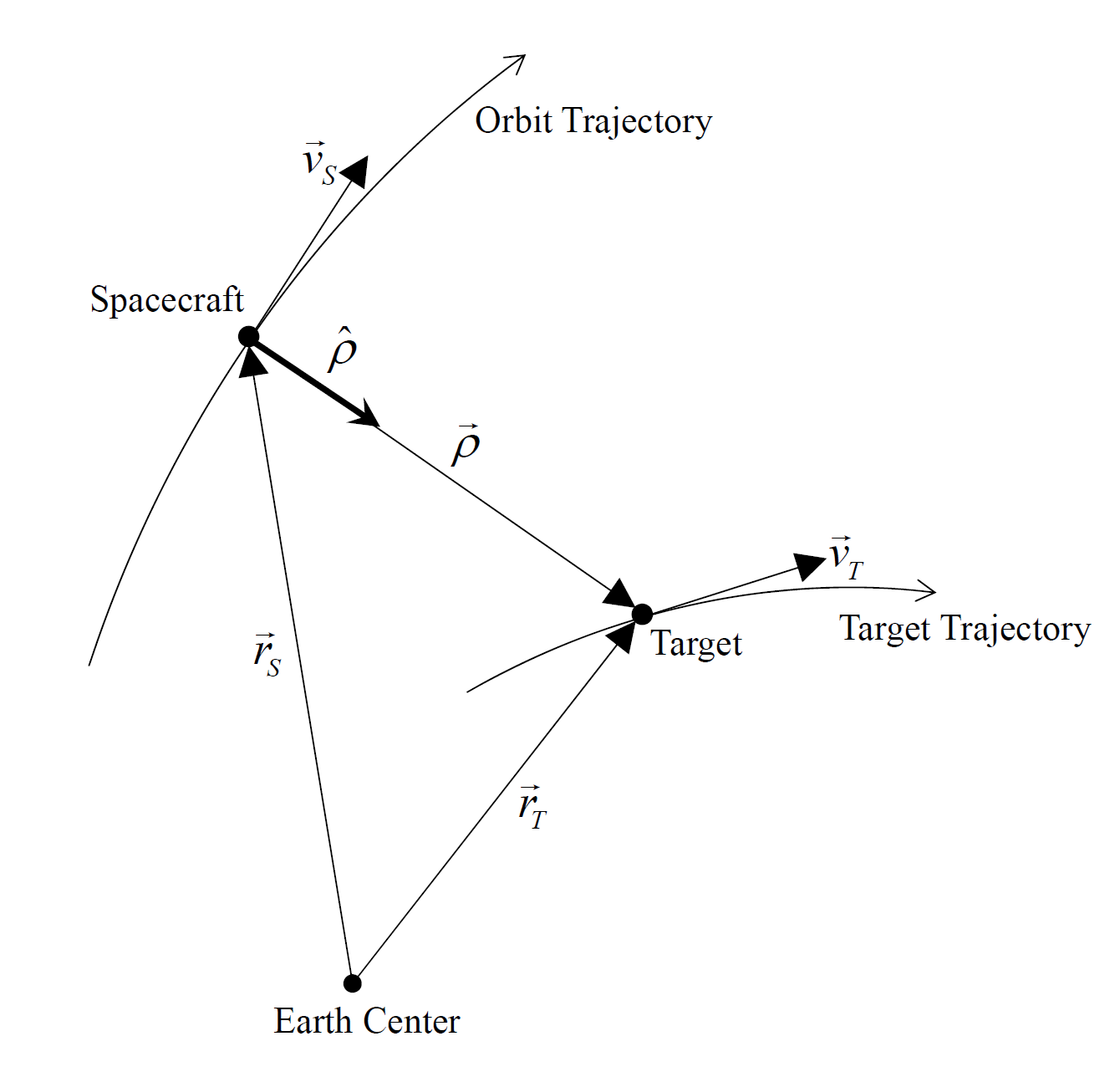}
\caption{Geometry of Staring Control.}
\label{fig_1}
\end{figure}

\begin{equation} \label{eqn_1} \begin{aligned}
\boldsymbol{\rho}(t) &= \textbf{\textit{r}}_T(t) - \textbf{\textit{r}}_S(t), &
\rho &= \lVert \textbf{\textit{r}}_{T} -\textbf{\textit{r}}_{S} \rVert, &
\hat{\boldsymbol{\rho}} &= \frac{\boldsymbol{\rho}}{\rho}
\end{aligned} \end{equation}
where $\textbf{\textit{r}}_T$ and $\textbf{\textit{r}}_S$ are the position vectors of the ground scan area and the satellite, respectively. Additional explanations related to $\textbf{\textit{r}}_T$ and $\textbf{\textit{r}}_S$ will be provided in the later section.

Let $\mathcal{D}$  be a command frame (the desired frame) that should be aligned with the body frame $\mathcal{B}$. Without loss of generality, the sensing axis of the payload is assumed to be aligned with the $z$-axis of the body frame $\mathcal{B}$.  Then, the starting condition to track the LOS vector by the desired $z$-axis vector $\hat{\textbf{\textit{z}}}_\mathcal{D}$
can be simply written as:
\begin{equation} \label{eq:z_D}
\hat{\textbf{\textit{z}}}_\mathcal{D} = \hat{\boldsymbol{\rho}}
\end{equation}
Note that, as mentioned earlier, Eq.~\eqref{eq:z_D} leaves the rotational degree of freedom along $\hat{\textbf{\textit{z}}}_\mathcal{D}$, unspecified. To fully define the desired attitude, either $\hat{\textbf{\textit{x}}}_\mathcal{D}$ or $\hat{\textbf{\textit{y}}}_\mathcal{D}$ must be specified. For convenience, $\hat{\textbf{\textit{x}}}_\mathcal{D}$ is here defined first, followed by $\hat{\textbf{\textit{y}}}_\mathcal{D}$, but the reverse sequence works as well. All the staring attitude commands can be generalized as follows:
\begin{equation}  \label{eqn_14}
\begin{aligned}
\hat{\textbf{\textit{z}}}_\mathcal{D} &= \boldsymbol{\rho}, & 
\hat{\textbf{\textit{x}}}_\mathcal{D} &= 
\frac{\hat{\textbf{\textit{z}}}_\mathcal{D} \times \textbf{\textit{k}} }{\lVert \hat{\textbf{\textit{z}}}_\mathcal{D} \times \textbf{\textit{k}} \rVert}, & 
\hat{\textbf{\textit{y}}}_\mathcal{D} &= \hat{\textbf{\textit{z}}}_\mathcal{D} \times \hat{\textbf{\textit{x}}}_\mathcal{D},
\end{aligned}
\end{equation}
where $\textbf{\textit{k}}$ represents the reference vector used for determining the desired attitude, and $\textbf{\textit{k}}$ must not $\hat{\textbf{\textit{z}}}_\mathcal{D} 	\parallel \textbf{\textit{k}}$. The reference vector can be designed to meet various criteria. Detailed explanations regarding the selection of $\textbf{\textit{k}}$ and the physical meaning of $\hat{\textbf{\textit{x}}}_\mathcal{D}$ for TDI camera will be provided in the following section. 

The analytic expression for the angular velocity of the general staring attitude is
\begin{equation} 
\begin{aligned}
\boldsymbol{\omega}_{\mathcal{D}} = 
\boldsymbol{\omega}_{z_\mathcal{D}^\perp} 
+ \boldsymbol{\omega}_{z_\mathcal{D}},\quad 
\boldsymbol{\omega}_{z_\mathcal{D}^\perp} 
=\frac{\boldsymbol{\rho} \times \dot{\boldsymbol{\rho}} }{\rho^2},
\quad 
\boldsymbol{\omega}_{z_\mathcal{D}} = \frac{\omega_{y_\mathcal{D}} k_{z_\mathcal{D}} - \dot{k}_{x_\mathcal{D}}}{k_{y_\mathcal{D}}} \hat{\textbf{\textit{z}}}_\mathcal{D}
\end{aligned}
\end{equation}
where $(\cdot)^\mathcal{D} = [(\cdot)_{x_\mathcal{D}},\ (\cdot)_{y_\mathcal{D}}, (\cdot)_{z_\mathcal{D}}]^\top$ is the expression of the vector $(\cdot)$ in $\mathcal{D}$. To be specific,
\begin{equation}
\begin{aligned}
    \boldsymbol{\omega}_{z_\mathcal{D}^\perp} &= \omega_{x_\mathcal{D}} \hat{\textbf{\textit{x}}}_\mathcal{D} + \omega_{y_\mathcal{D}} \hat{\textbf{\textit{y}}}_\mathcal{D} \\
    \boldsymbol{\omega}_{\mathcal{D}} &= \omega_{x_\mathcal{D}} \hat{\textbf{\textit{x}}}_\mathcal{D} + \omega_{y_\mathcal{D}} \hat{\textbf{\textit{y}}}_\mathcal{D} + \omega_{z_\mathcal{D}} \hat{\textbf{\textit{z}}}_\mathcal{D}
\end{aligned}
\implies
\begin{aligned}
    \boldsymbol{\omega}_{z_\mathcal{D}^\perp}^{\mathcal{D}} &= [\omega_{x_\mathcal{D}},\ \omega_{y_\mathcal{D}},\ 0]^\top \\
    \boldsymbol{\omega}_{\mathcal{D}}^{\mathcal{D}} &= [\omega_{x_\mathcal{D}},\ \omega_{y_\mathcal{D}},\ \omega_{z_\mathcal{D}}]^\top \\
\end{aligned}
\end{equation}
Likewise, the angular acceleration becomes
\begin{equation} 
\begin{aligned}
\boldsymbol{\alpha}_{\mathcal{D}} &= 
\boldsymbol{\alpha}_{z_\mathcal{D}^\perp} 
+ \boldsymbol{\alpha}_{z_\mathcal{D}} \\
\boldsymbol{\alpha}_{z_\mathcal{D}^\perp} 
&= \frac{\boldsymbol{\rho} \times \ddot{\boldsymbol{\rho}}}{\rho^2}
- 2\frac{ \boldsymbol{\rho} \cdot \dot{\boldsymbol{\rho}}}{\rho^2} \boldsymbol{\omega}_{z_\mathcal{D}^\perp} + 
\boldsymbol{\omega}_{z_\mathcal{D}^\perp} \times \boldsymbol{\omega}_{z_\mathcal{D}}\\
\boldsymbol{\alpha}_{z_\mathcal{D}} &=\frac{ \alpha_{y_\mathcal{D}} k_{z_\mathcal{D}}
- \omega_{x_\mathcal{D}} \omega_{z_\mathcal{D}} k_{z_\mathcal{D}}
- \omega_{x_\mathcal{D}} \omega_{y_\mathcal{D}} k_{y_\mathcal{D}}
+ 2\omega_{y_\mathcal{D}}\dot{k}_{z_\mathcal{D}}
- 2\omega_{z_\mathcal{D}} \dot{k}_{y_\mathcal{D}}
- \ddot{k}_{x_\mathcal{D}}}{k_{y_\mathcal{D}}}\hat{\textbf{\textit{z}}}_\mathcal{D}
\end{aligned}
\end{equation}
Note that the profile of $\textbf{\textit{r}}_S$, $\textbf{\textit{r}}_T$, $\textbf{\textit{k}}$ and their first and second derivatives are necessary to compute the analytic attitude profile\cite{han2022analytical}.

\subsection{General Strip Attitude for TDI camera}
\begin{figure}[h]
\centering
\includegraphics[width=0.4\textwidth]{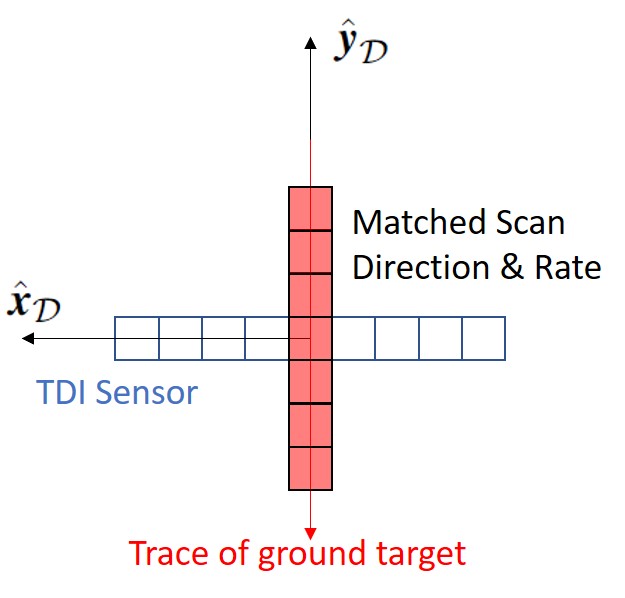}
\includegraphics[width=0.48\textwidth]{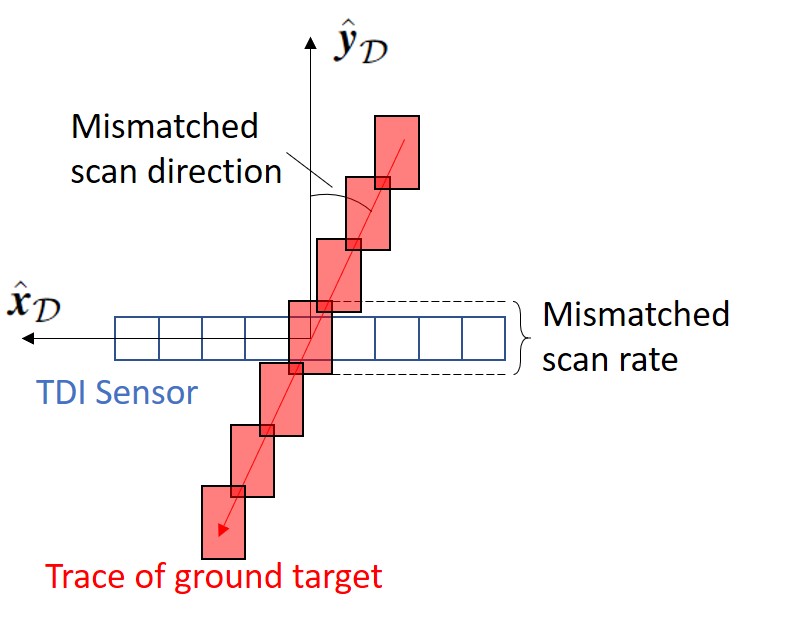}
\caption{Illustration of (a) matched (b) mismatched scan direction and scan rate errors}
\label{fig_2}
\end{figure}

To obtain a high-quality image using a TDI camera, it is necessary to synchronize the relative motion between the camera and the object. This motion synchronization can be described by two factors, as shown in Figure~\ref{fig_2}: the line scan direction, which determines whether the direction of relative motion aligns with the sensor's scan direction, and the line scan rate, which determines whether the frequency of each exposure matches the speed of relative motion. Without loss of generality, we will assume the scan direction of the sensor is aligned with
$\hat{\textbf{\textit{y}}}_\mathcal{B}$, meaning that $\hat{\textbf{\textit{y}}}_\mathcal{D}$ should be designed to match the scan direction.

Both scan direction and scan rate are related to the relative motion of the current ground target (i.e., the instantaneous ground point on the boresight axis) in the sensor frame\cite{yoon2009analytical, han2017attitude}. To analyze the relative motion, we first define the line of sight vector between the satellite and the instantaneous ground target at $t_i$, $\textbf{\textit{r}}_T(t_i)$ as follows:
\begin{equation}
    \boldsymbol{\rho}_i(t) = \textbf{\textit{r}}_T(t_i) - \textbf{\textit{r}}_S(t)
\end{equation}
Then the relative motion between the satellite and the ground target at $t_i$ in the Earth-centered inertial (ECI) frame $\mathcal{I}$ becomes:
\begin{equation}
    {}^\mathcal{I} \dot{\boldsymbol{\rho}}_i(t) = \boldsymbol{\omega}_\mathcal{F/I}\times\textbf{\textit{r}}_T(t_i) - \textbf{\textit{v}}_S(t)
\end{equation}
since the target is fixed in $\mathcal{F}$, which is the Earth-centered, Earth-fixed (ECEF) frame, and $\boldsymbol{\omega}_\mathcal{F/I}$ is the rotational vector of the Earth. Note that the left superscript of the vector derivative might be omitted if the inertial frame is used.

The relative motion with respect to $\mathcal{D}$ then becomes:
\begin{equation}
\begin{aligned}
    {}^\mathcal{D} \dot{\boldsymbol{\rho}}_i(t) &= {}^\mathcal{I} \dot{\boldsymbol{\rho}}_i(t) + \boldsymbol{\omega}_\mathcal{D/I} \times \boldsymbol{\rho}_i(t) \\
    &= \boldsymbol{\omega}_\mathcal{F/I}\times\textbf{\textit{r}}_T(t_i) - \textbf{\textit{v}}_S(t)
    + \left( \boldsymbol{\omega}_{z_\mathcal{D}^\perp } (t)
+ \boldsymbol{\omega}_{z_\mathcal{D}}(t) \right) \times \boldsymbol{\rho}_i(t)
\end{aligned}    
\end{equation}
Using the fact that $\boldsymbol{\rho}_i(t_i) = \boldsymbol{\rho}(t_i)$ and $\boldsymbol{\omega}_{z_\mathcal{D}}(t) \parallel \hat{\textbf{\textit{z}}}_\mathcal{D}$,  ${}^D \dot{\boldsymbol{\rho}}_i(t)$ at $t=t_i$ becomes:
\begin{equation}
\begin{aligned}
    {}^\mathcal{D} \dot{\boldsymbol{\rho}}_i(t_i) 
    &= \boldsymbol{\omega}_\mathcal{F/I}\times\textbf{\textit{r}}_T(t_i) - \textbf{\textit{v}}_S(t_i)
    + \frac{\boldsymbol{\rho}(t_i) \times \dot{\boldsymbol{\rho}}(t_i) }{\rho^2(t_i)} \times \boldsymbol{\rho}_i(t_i) \\
    &= 
    \left\{\dot{\boldsymbol{\rho}}(t_i) - {}^\mathcal{F} \dot{\textbf{\textit{r}}}_T(t_i) \right\}+ 
    \left\{ -\dot{\boldsymbol{\rho}}(t_i) + 
    \left( \hat{\boldsymbol{\rho}}(t_i) \cdot \dot{\boldsymbol{\rho}}(t_i)\right)\hat{\boldsymbol{\rho}}(t_i) \right\} \\
    &= 
    - {}^\mathcal{F} \dot{\textbf{\textit{r}}}_T(t_i) + \left( \hat{\textbf{\textit{z}}}_\mathcal{D}(t_i) \cdot \dot{\boldsymbol{\rho}}(t_i)\right)\hat{\textbf{\textit{z}}}_\mathcal{D} (t_i)
\end{aligned}    
\end{equation}
Therefore, define $\textbf{\textit{k}}$ as
\begin{equation}
    \textbf{\textit{k}} = - {}^\mathcal{F} \dot{\textbf{\textit{r}}}_T = -\textbf{\textit{v}}_T + \boldsymbol{\omega}_\mathcal{F/I} \times \textbf{\textit{r}}_T
\end{equation}
will align the relative motion with $\hat{\textbf{\textit{y}}}_\mathcal{D}$, making a zero drift angle. Its first and second inertial derivatives are:
\begin{equation}
\begin{aligned}
    {}^\mathcal{I} \dot{\textbf{\textit{k}}} &= - \left(\textbf{\textit{a}}_T + \boldsymbol{\omega}_\mathcal{F/I}\times\textbf{\textit{v}}_T \right) \\
    {}^\mathcal{I} \ddot{\textbf{\textit{k}}} &= - \left( \dot{\textbf{\textit{a}}}_T +  \boldsymbol{\omega}_\mathcal{F/I}\times\textbf{\textit{a}}_T \right)
\end{aligned}
\end{equation}
Note that the time derivative of target acceleration, i.e., jerk, is required to compute the second derivative.

On the other hand, projecting the ${}^\mathcal{D} \dot{\boldsymbol{\rho}}_i(t_i)$ onto the sensor plate $\hat{\textbf{\textit{x}}}_\mathcal{D}-\hat{\textbf{\textit{y}}}_\mathcal{D}$ gives:
\begin{equation}
\begin{aligned}
    {}^\mathcal{D} \dot{\boldsymbol{\rho}}_i - \left({}^\mathcal{D} \dot{\boldsymbol{\rho}}_i \cdot \hat{\textbf{\textit{z}}}_\mathcal{D}\right) \hat{\textbf{\textit{z}}}_\mathcal{D}
    &= - {}^\mathcal{F} \dot{\textbf{\textit{r}}}_T + \left( \hat{\textbf{\textit{z}}}_\mathcal{D}\cdot \dot{\boldsymbol{\rho}}\right)\hat{\textbf{\textit{z}}}_\mathcal{D}  +
    \left({}^\mathcal{F} \dot{\textbf{\textit{r}}}_T \cdot \hat{\textbf{\textit{z}}}_\mathcal{D}\right) \hat{\textbf{\textit{z}}}_\mathcal{D} - \left( \hat{\textbf{\textit{z}}}_\mathcal{D}\cdot \dot{\boldsymbol{\rho}}\right)\hat{\textbf{\textit{z}}}_\mathcal{D} \\
    &= - {}^\mathcal{F} \dot{\textbf{\textit{r}}}_T + 
    \left({}^\mathcal{F} \dot{\textbf{\textit{r}}}_T \cdot \hat{\textbf{\textit{z}}}_\mathcal{D}\right) \hat{\textbf{\textit{z}}}_\mathcal{D}  \\
\end{aligned}
\end{equation}
and the norm of the projected vector represents the speed of relative motion on the sensor plane as:
\begin{equation}
\begin{aligned}
    v_\text{LOS} &= \lVert {}^D \dot{\boldsymbol{\rho}}_i - \left({}^D \dot{\boldsymbol{\rho}}_i \cdot \hat{\textbf{\textit{z}}}_\mathcal{D}\right)\hat{\textbf{\textit{z}}}_\mathcal{D} \rVert \\
    &=  \sqrt{
    \lVert {}^F \dot{\textbf{\textit{r}}}_T \rVert^2 + 
    \left({}^F \dot{\textbf{\textit{r}}}_T \cdot \hat{\textbf{\textit{z}}}_\mathcal{D}\right)^2 - 
    2\left( {}^F \dot{\textbf{\textit{r}}}_T \cdot \hat{\textbf{\textit{z}}}_\mathcal{D}\right)^2 }\\
    &= \lVert {}^F \dot{\textbf{\textit{r}}}_T \rVert \sin \psi
\end{aligned}
\end{equation}
where $\psi$ is the angle between ${}^F \dot{\textbf{\textit{r}}}_T$ and $\hat{\textbf{\textit{z}}}_\mathcal{D}$. The projected speed can be converted to the line scan speed and exposure frequency for given camera parameters as:
\begin{equation}
    v_\text{CCD}(t) = \frac{p_\text{CCD}}{\rho(t)}v_\text{LOS}(t),\quad f_\text{CCD}(t) = \frac{v_\text{CCD}(t)}{d_\text{CCD}}
\end{equation}
where $p_\text{CCD}$ is the focal length of camera and $d_\text{CCD}$ is the pixel size.

\begin{figure}[h]
\centering
\includegraphics[width=0.35\textwidth]{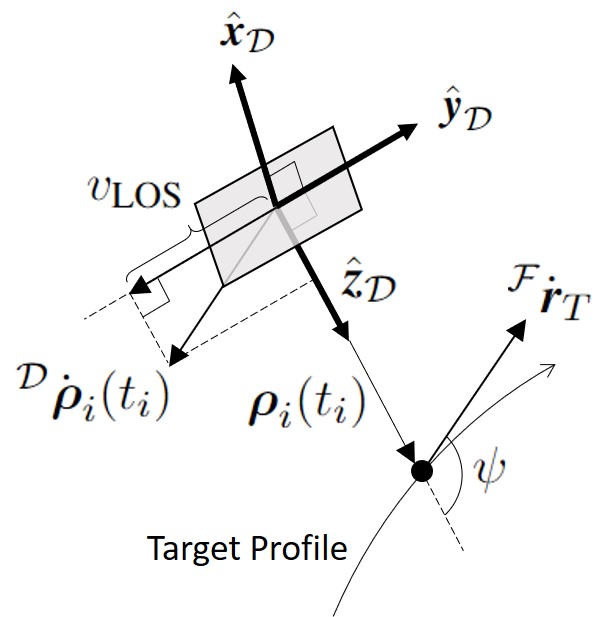}
\includegraphics[width=0.35\textwidth]{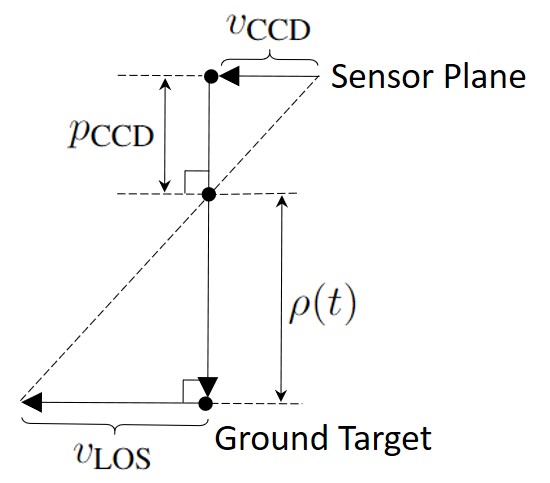}
\caption{(a) Correctly aligned frame $\mathcal{D}$ and (b) relationship between $v_\text{LOS}$ and $v_\text{ccd}$}
\label{fig_3}
\end{figure}

In summary, for a given ground profile $\textbf{\textit{r}}_T(t)$, the proposed $\textbf{\textit{k}}$ and  $v_\text{CCD}$ will give a zero drift angle attitude with a synchronized line rate. The additional goal of this research is to find the optimal $\textbf{\textit{r}}_T(t)$ if we have freedom to select it.

\subsection{Target Position Model}
Let's assume the target area of interest is determined and known by either the ground or space system. In practice, a great circle of spherical earth model, or loxodrome/geodesic with an ellipsoidal earth model, or the general curve fitted parameterized curve, or a combination of them are used to model a ground target. In all cases, they can be represented as:
\begin{equation}
    \textbf{\textit{r}}_T(t) =\textbf{\textit{r}}_T(s(t), \theta)
\end{equation}
where $s$ is parameter profile and and $\theta$ are the constant parameters defining the curve. 

\begin{figure}[h]
\centering
\includegraphics[width=0.35\textwidth]{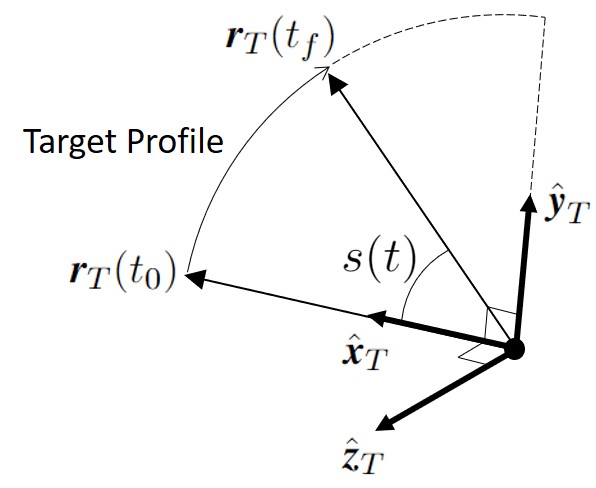}
\caption{Target profile as a great circle of spherical earth model}
\label{fig_4}
\end{figure}

For simplicity, this paper uses the great circle of a spherical Earth model. Let's assume the starting and ending positions of the strip are given as $\textbf{\textit{r}}_T(t_0)$ and $\textbf{\textit{r}}_T(t_f)$, respectively. Then, the great circle can be parameterized as:
\begin{equation}
    \textbf{\textit{r}}_T(s(t)) = R_E \cos s(t)\ \hat{\textbf{\textit{x}}}_T + R_E \sin s(t)\ \hat{\textbf{\textit{y}}}_T
\end{equation}
where $R_E$ is the radius of spherical earth, $\hat{\textbf{\textit{x}}}_T$ is the radial direction for the starting position of STRIP, $\hat{\textbf{\textit{y}}}_T$ is the along-track direction of STRIP, and $s(t)$ is the traveling angle from $\hat{\textbf{\textit{x}}}_T$. Referring to Figure.~\ref{fig_4}, the $\hat{\textbf{\textit{x}}}_T$, $\hat{\textbf{\textit{y}}}_T$, and $\hat{\textbf{\textit{z}}}_T$ are defined as:
\begin{equation}
    \hat{\textbf{\textit{x}}}_T = \frac{\textbf{\textit{r}}_T(t_0)}{r_T(t_0)}, \quad
    \hat{\textbf{\textit{z}}}_T = \frac{\textbf{\textit{r}}_T(t_0) \times \textbf{\textit{r}}_T(t_f)}{\lVert \textbf{\textit{r}}_T(t_0) \times \textbf{\textit{r}}_T(t_f) \rVert}, \quad
    \hat{\textbf{\textit{y}}}_T = \hat{\textbf{\textit{z}}}_T \times \hat{\textbf{\textit{x}}}_T
\end{equation}
and the range of $s$ is $s(t) \in [0, \cos^{-1}(\textbf{\textit{r}}_T(t_0) \cdot \textbf{\textit{r}}_T(t_f)/R_E^2)]$. 

The velocity and acceleration of the ground target with respect to $\mathcal{F}$ can be computed as:
\begin{equation}
    {}^\mathcal{F} \dot{\textbf{\textit{r}}}_T = 
    R_E \dot{s}\left(  - \sin s\ \hat{\textbf{\textit{x}}}_T + \cos s\ \hat{\textbf{\textit{y}}}_T \right)
\end{equation}
and
\begin{equation}
    {}^\mathcal{F} \ddot{\textbf{\textit{r}}}_T = 
    R_E \ddot{s} \left(  - \sin s\ \hat{\textbf{\textit{x}}}_T + \cos s\ \hat{\textbf{\textit{y}}}_T \right)
    - R_E \dot{s}^2\left(  \cos s\ \hat{\textbf{\textit{x}}}_T + \sin s\ \hat{\textbf{\textit{y}}}_T \right)
\end{equation}
One can compute the derivatives with respect to $\mathcal{I}$ as:
\begin{equation}
\begin{aligned}
    \textbf{\textit{v}}_T &= {}^\mathcal{I} \dot{\textbf{\textit{r}}}_T =  
    {}^\mathcal{F} \dot{\textbf{\textit{r}}}_T + \boldsymbol{\omega}_\mathcal{F/I} \times \textbf{\textit{r}}_T \\
    \textbf{\textit{a}}_T &= {}^I \dot{\textbf{\textit{v}}}_T
    = {}^\mathcal{F} \ddot{\textbf{\textit{r}}}_T + 2\boldsymbol{\omega}_\mathcal{F/I} \times {}^\mathcal{F} \dot{\textbf{\textit{r}}}_T
    + \boldsymbol{\omega}_\mathcal{F/I} \times \left( \boldsymbol{\omega}_\mathcal{F/I} \times \textbf{\textit{r}}_T \right)
\end{aligned}
\end{equation}

\subsection{Satellite Orbit Model}
As long as reasonable time-stamped position and velocity data are available, there is no restriction on the satellite orbit model for the previously explained attitude equation. When deriving the angular acceleration, the acceleration of the satellite is required, but simply using a two-body gravitational model (or including up to the $J_2$ term) is accurate enough for this application.

Depending on the architecture of operation management, orbit prediction may or may not be required by the onboard computer. Based on previous work\cite{kim2019sensitivity}, in-orbit accuracy of 5 m for position and 0.6 cm/s for velocity can be achieved with the aid of a high-fidelity gravitational model and GPS data. With this initial orbital accuracy, the accuracy of orbit prediction can be maintained using simple numerical integration with sparse time steps if the prediction horizon is about a couple of minutes, allowing in-orbit attitude profile generation to be feasible. If the initial and terminal times for the imaging operation are fixed, only a single orbit prediction is required.

\section{Constrained Optimal Control Problem and DDP}

\subsection{Unconstrained Optimal Control Problem and the HJB equation}
The general trajectory optimization with terminal state constraints and fixed terminal time can be formulated as follows:
\begin{equation}
\begin{aligned}
    \min_{u} \ &J = \phi(x(t_f), t_f) + \int_{t_0}^{t_f} L(x(t),u(t),t)\ dt \\
    \textrm{s.t.} \ &\dot x = f(x(t), u(t),t), \ x(t_0) = x_0, \ \psi(x(t_f), t_f) = 0 \\
\end{aligned}
\end{equation}
where $x\in \mathbb{R}^n$ is the states, $u\in \mathbb{R}^m$ is the control, $\phi:\mathbb{R}^n\times\mathbb{R}\to\mathbb{R}$ is the terminal cost, 
$L:\mathbb{R}^n\times\mathbb{R}^m\times\mathbb{R}\to\mathbb{R}$ is the learning cost, 
$f:\mathbb{R}^n\times\mathbb{R}^m\times\mathbb{R}\to\mathbb{R}^n$ is the state dynamics, and
$\psi:\mathbb{R}^n\times\mathbb{R}\to\mathbb{R}^d$ is the terminal constraints.

Let $V(x,t) \triangleq \min J(x,t,u)$ and assume there exist twice differentiable $V$ satisfying the following HJB equation:
\begin{equation}
    \begin{aligned}
        - V_t (x,t) = \min_{u} \left\{ L(x,u,t) + V_x(x,t)\cdot f(x,u,t)\right\}
    \end{aligned}
\end{equation}
then for $\theta(x,t)$ given by:
\begin{equation}
    \begin{aligned}
        \theta (x,t) = \argmin_{u} \left\{ L(x,u,t) + V_x(x,t)\cdot f(x,u,t)\right\}
    \end{aligned}
\end{equation}
is the optimal control as:
\begin{equation}
    u^\ast(t) = \theta(x^\ast,t)
\end{equation}
where $x^\ast$ satisfies
\begin{equation}
    \dot x^\ast = f(x^\ast, u^\ast, t)
\end{equation}
Note that in order to handle the terminal constraints $\psi$, it is proposed that the Lagrange multiplier $\nu\in \mathbb{R}^d$ is employed to adjoin the terminal constraint into cost function \cite{b1}.
\begin{equation}
    \Tilde{\phi} (x(t_f), t_f,\nu) = \phi (x(t_f), t_f) + \nu^\top \psi(x(t_f), t_f)
\end{equation}
then the terminal condition for cost to go function $V$ becomes:
\begin{equation}
    V(x(t_f),t_f, \nu) = \Tilde{\phi} (x(t_f), t_f,\nu)
\end{equation}

\subsection{Review of Differential Dynamic Programming}
The objective of DDP is to numerically compute $u^\ast(t)$ using the second-order approximation of the HJB equation. There are various ways to handle terminal and path constraints; this paper handles the terminal constraint as one of the DDP parameters and the path constraint by adding a penalty cost representing the violation of constraints.

\subsubsection{DDP Updated Algorithm with Terminal Constraint} 
To approximate the HJB equations, all functions are expanded up to the second order with respect to the reference trajectory ($\bar x,\ \bar u$) and the reference parameters ($\bar\nu$). With the reference values, the variations of state, control, and parameters are defined as follows:
\begin{equation}
    \delta x \triangleq x^\ast - \bar{x},\quad
    \delta u \triangleq u^\ast - \bar{u},\quad
    \delta \nu \triangleq \nu^\ast - \bar{\nu}
\end{equation}
With this definition, the linearization of the dynamics becomes:
\begin{equation} \label{eq:forward}
\begin{aligned}
    \frac{d}{dt} ({\bar{x}+\delta x}) = f({\bar{x}+\delta x}, {\bar{u}+\delta u}, t) \implies
     \dot{\delta x} \approx \bar f_x \delta x + \bar f_u \delta u 
\end{aligned}
\end{equation}
and the second-order approximation of the cost-to-go function becomes:
\begin{equation}
\begin{aligned}
    V(x^\ast, \nu^\ast, t)
    &= V(\bar x + \delta x, \bar \nu + \delta \nu, t) \\
     &\approx \bar V + \bar V_x \cdot \delta x + \bar V_\nu \cdot \delta \nu  
     + \frac{1}{2} 
     \begin{bmatrix} \delta x \\ \delta \nu \end{bmatrix}^\top
     \begin{bmatrix} \bar V_{xx} & \bar V_{x\nu} \\ \bar V_{\nu x} & \bar V_{\nu \nu} \end{bmatrix}
     \begin{bmatrix} \delta x \\ \delta \nu \end{bmatrix}
\end{aligned}
\end{equation}
Then, the update equation for each term can be derived as follows\cite{sun2014continuous}:
\begin{equation} \label{eq:backward}
\begin{aligned}
    -\dot {\bar{V}} &= \bar{L} + \bar Q_u \beta_u + \frac{1}{2} \beta_u^\top \bar Q_{uu} \beta_u \\
    -\dot {\bar{V}}_x &= \bar{Q}_x + \bar{Q}_u \beta_x + \beta_u \bar Q_{ux} +  \beta_u^\top \bar Q_{uu} \beta_x \\
    -\dot {\bar{V}}_\nu &= \bar{Q}_u \beta_\nu + \beta_u^\top \bar Q_{uu} \beta_\nu + \beta_u^\top \bar Q_{u\nu} \\
    -\dot {\bar{V}}_{xx} &= \bar{Q}_{xx} + 2\bar{Q}_{xu}\beta_x + \beta_x^\top \bar{Q}_{uu} \beta_x \\
    -\dot {\bar{V}}_{x\nu} &= \bar{Q}_{x\nu} + \bar{Q}_{xu}\beta_\nu + \beta_u^\top \bar Q_{uu} \beta_\nu \\
    -\dot {\bar{V}}_{\nu \nu} &= \beta_\nu^\top \bar{Q}_{uu} \beta_\nu + 2 \bar{Q}_{\nu u} \beta_\nu\\    
\end{aligned}
\end{equation}
where the $Q$ matrices are defined as:
\begin{equation} \label{eq:q_matrix}
\begin{aligned}
    &\bar Q_x \triangleq \bar L_x + \bar V_x^\top \bar f_x,\quad
    \bar Q_u \triangleq \bar L_u + \bar V_x^\top \bar f_u,\quad
    \bar Q_{xx} \triangleq \bar L_{xx} + 2\bar{V}_{xx}\bar f_x \\
    &\bar Q_{xu} \triangleq \bar L_{xu} + \bar{V}_{xx}\bar f_u,\quad 
    \bar Q_{\nu x} \triangleq \bar V_{\nu x} \bar f_x,\quad
    \bar Q_{uu} \triangleq \bar L_{uu},\quad 
    \bar Q_{\nu u} \triangleq \bar V_{\nu x} \bar f_u \\ 
\end{aligned}
\end{equation}
and the control gains $\beta$ are defined as:
\begin{equation}
    \bar \beta_u \triangleq -\bar{Q}_{uu}^{-1}\bar Q_u,\quad
    \bar \beta_x \triangleq -\bar{Q}_{uu}^{-1}\bar Q_{ux},\quad
    \bar \beta_\nu \triangleq -\bar{Q}_{uu}^{-1}\bar Q_{u \nu} \\
\end{equation}
with the terminal conditions as:
\begin{equation} \label{eq:terminal_condition}
\begin{aligned}
    &\bar V( t_f) = \bar \phi + \bar \nu^\top \bar \psi,\quad 
    \bar V_x( t_f) =\bar \phi_x + \bar \nu^\top \bar \psi_x,\quad 
    \bar V_\nu( t_f) = \bar \psi^\top\\
    &\bar V_{xx}( t_f) =  \bar \phi_{xx} + \bar\nu^\top \bar \psi_{xx},\quad 
    \bar V_{x\nu}( t_f) = \bar \psi_x ^\top,\quad 
    \bar V_{\nu\nu}( t_f)=  0
\end{aligned}
\end{equation}

After completing the backward integration of Eq.~\eqref{eq:backward}, one can compute the variation of $V$ at $t_0$ as
\begin{equation}
\begin{aligned}
    V(x^\ast, \nu^\ast, t_0) \approx \bar V(t_0) + \bar V_\nu(t_0)\delta \nu + \frac{1}{2} \delta \nu^\top\bar V_{\nu \nu}(t_0) \delta \nu
\end{aligned}
\end{equation}
since $\delta x(t_0)=0$. Therefore, the correction for $\nu$ is set to minimize the $V$ as:
\begin{equation} \label{eq:nu_correction}
    \delta \nu ^\ast = -\bar V_{\nu \nu}^{-1}(t_0)\bar V_\nu(t_0) \implies \bar \nu^+ = \bar \nu^- + k_\nu \delta \nu ^\ast
\end{equation}
Likewise, the correction terms for $\bar u$ are computed iteratively as:
\begin{equation} \label{eq:u_correction}
    \delta u^\ast 
    = \beta_u +  \beta_x \delta x + \beta_\nu \delta \nu^\ast \implies
    \bar u^+ = \bar u^- + k_u \delta u^\ast
\end{equation}
with the Eq.~\eqref{eq:forward}. The amount of correction is controlled by the gains 
$k_\nu$, $k_{u} \in [0,1)$.

\subsubsection{Augmented Lagrangian Method for Path Constraint}
This paper adopts the augmented Lagrangian method to handle the path constraints, and a detailed explanation can be found in the referenced work\cite{plancher2017constrained, howell2019altro}. We will briefly explain the overall strategy of the method.

Let $g(x,u)\leq 0$ be the inequality path constraints and $h(x,u)=0$ be the equality path constraints. Then, augmented learning cost is defined as\cite{toussaint2014novel}:
\begin{equation}
    L_A = L + \sum_{i} (\lambda_i g_i + \mathbf{1}_i \mu_i g_i^2) + \sum_{i} (\eta_i h_i + \kappa_i h_i^2)
\end{equation}
where $\lambda_i$ and $\eta_i$ are Lagrange multipliers, and ${\mu}_i$ and $\kappa_i$ are penalty weights. The indicator function for inequality constraint is defined as:
\begin{equation}
    \mathbf{1}_i  = \begin{cases}
        1 & \mbox{if } g_i \geq 0\ \vee \ \lambda_i > 0 \\
        0 & \mbox{else }
    \end{cases}
\end{equation}
The Q matrices in Eq.~\eqref{eq:q_matrix} are updated with the new learning cost, and then the previously explained DDP method can be used.

After conducting the DDP logic by holding ($\lambda$, $\mu$) and ($\eta$, $\kappa$) constant, the Lagrange multipliers are updated as follows:
\begin{equation} \label{eq:update_lagrange}
    \lambda_i^+ = \max(0, \lambda_i^- + \mu_i g_i),\quad
    \eta_i^+ = \eta_i^- + \kappa_i h_i
\end{equation}
and the penalty weights are increased monotonically with a scaling factor $\phi > 1$ as:
\begin{equation} \label{eq:update_penalty}
    \mu_i^+ = \gamma \mu_i^-,\quad \kappa_i^+ = \gamma \kappa_i^-
\end{equation}

The overall procedure of CDDP is summarized in the following tables.
\begin{algorithm}
    \caption{Constrained DDP Algorithm}
  \begin{algorithmic}[1]
    \REQUIRE Optimal Control Problem $f$, $\psi$, $\phi$, $L$, $g$, $h$ 
    \INPUT Initial Guess ($\bar x$, $\bar{u}$, $\bar{\nu}$), Initial Lagrangian and Penalty ($\lambda$, $\mu$, $\eta$, $\kappa$), Update Parameters ($k_u$, $k_\nu$, $\gamma$)
    \OUTPUT Optimal Trajectory ($x^\ast$, $u^\ast$, $\nu^\ast$)
    \WHILE{$ \lvert\delta V \rvert > \epsilon_V $ or $\lvert h \rvert >\epsilon_h$ or $ g >\epsilon_g$}
      \STATE Generate the reference trajectory $\bar{x}(t)$ using $\bar{u}(t)$\
      \STATE Integrate the cost-to-go function and its derivatives backward using Eqs.~\eqref{eq:backward}-\eqref{eq:terminal_condition}
      \STATE Update $\bar \nu$ using Eq.~\eqref{eq:nu_correction} and $\bar u$ using Eq.~\eqref{eq:u_correction}
      \STATE Update Penalty Parameters using Eqs.~\eqref{eq:update_lagrange}-\eqref{eq:update_penalty}
    \ENDWHILE
  \end{algorithmic}
\end{algorithm}

\subsection{Optimal Strip Attitude Command Generation}
If the ground scan profile is given, then the 3-axis attitude of the EOS is determined for the entire imaging duration. However, if the TDI camera can adjust the $f_\text{CCD}(t)$ during the imaging operation, we can control how fast the EOS should scan the ground target profile. We aim to design $f_\text{CCD}(t)$ such that it minimizes the angular velocity during the imaging operation. Additionally, we also consider the case where there is a hardware limit for $f_\text{CCD}$, which can be modeled as $f_\text{LB} \leq f_\text{CCD} \leq f_\text{UB}$.

In this paper, we have considered the following two OCP, which are
\begin{equation}
\begin{aligned}
    \min_{u(\cdot)} \quad & \int_{t_0}^{t_f} \boldsymbol{\omega}(t, s(t), u(t)) \cdot \boldsymbol{\omega}(t, s(t), u(t))\ dt\\
    \textrm{subject to} \quad & \dot s = u, \quad s(t_0) = 0,\quad s(t_f) = \cos^{-1}(\textbf{\textit{r}}_T(t_0) \cdot \textbf{\textit{r}}_T(t_f)/R_E^2)\\
    \quad & f_\text{LB} \leq f_\text{CCD}(t) \leq f_\text{UB}\\
\end{aligned} 
\end{equation}
and
\begin{equation}
\begin{aligned}
    \underset{u(\cdot),\ t\in [t_0,t_f]}{\text{min\ max}}\quad & \lVert \boldsymbol{\omega}(t, s(t), u(t)) \rVert \\
    \textrm{subject to} \quad & \dot s = u, \quad s(t_0) = 0,\quad s(t_f) = \cos^{-1}(\textbf{\textit{r}}_T(t_0) \cdot \textbf{\textit{r}}_T(t_f)/R_E^2)\\
    \quad & f_\text{LB} \leq f_\text{CCD}(t) \leq f_\text{UB}\\
\end{aligned} 
\end{equation}
The first OCP minimizes the integral of the rate squared, whereas the second OCP minimizes the maximum angular velocity. Note that both performance indices are directly related to the quality of images. Due to the complex expression for $\boldsymbol{\omega}$ and the generalized structure of both satellite and ground scan position models, finding an analytic solution is intractable unless simplified models are used. Therefore, we propose to solve the OCP using the DDP method.

\textbf{Remark:} Note that the current formulation $\dot s = u$ with zero-order hold discretization cannot reflect the learning cost at the last time step. This issue will be addressed in future work.

From the previous section, it is necessary to ensure the positive definiteness of $Q_{uu}$, which mainly depends on $L_{uu}$. We can roughly analyze the $L_{uu}$ by following approximation:
\begin{equation}
\begin{aligned}
    L &\approx \boldsymbol{\omega}_{z_\mathcal{D}^\perp} \cdot \boldsymbol{\omega}_{z_\mathcal{D}^\perp} 
    = \frac{1}{\rho^4} \left( (\boldsymbol{\rho} \times \dot{\boldsymbol{\rho}}) \cdot (\boldsymbol{\rho} \times \dot{\boldsymbol{\rho}}) \right)
    = \frac{1}{\rho^2} \left( \dot{\rho}^2  - (\hat{\boldsymbol{\rho}} \cdot \dot{\boldsymbol{\rho}})^2 \right) 
\end{aligned}
\end{equation}
and this approximation is reasonable for most scenarios of interest. Re-expressing the target velocity as $\textbf{\textit{v}}_T = u \boldsymbol{\nu}_T$ gives
 \begin{equation}
     L \approx \frac{1}{\rho^2} \left( 
     \left( \nu_T^2 - (\hat{\boldsymbol{\rho}}\cdot \boldsymbol{\nu}_T)^2  \right) u^2 + \cdots \right) \implies L_{uu} \approx \frac{\nu_T^2 - (\hat{\boldsymbol{\rho}}\cdot \boldsymbol{\nu}_T)^2 }{\rho^2} \geq 0 
 \end{equation}
Note that $L_{uu}$ becomes zeros when $\boldsymbol{\rho} \parallel \boldsymbol{\nu}_T$ meaning that the boresight axis is parallel to the target profile direction. This only happens when the target is located at the horizon as seen from the satellite, which is not a case of interest. Although this is an approximated analysis, we can expect reasonable behavior of $L_{uu}$ for most scenarios.

To solve the second OCP, we first approximate the maximum angular velocity using the softmax function as:
\begin{equation}
    \max_{t \in [t_0, t_f]} \omega(t) \propto \ln \left( \int_{t_0}^{t_f} \exp(\omega(t)) \ dt \right)
    \propto \ln \ln \left( \int_{t_0}^{t_f} \exp\left(\exp(\omega(t))\right) \ dt \right) \approx \cdots
\end{equation}
Note that the more exponential terms used, the more accurate the solution. Due to the strict monotonicity of the logarithm function, minimizing the integral will yield an equivalent result, allowing us to apply the previously explained DDP algorithm to find the optimal trajectory. To effectively find the solution while avoiding numerical issues, appropriate constants $N$ and $M$ are used and updated every iteration based on the $\omega(t)$ profile of the previous iteration. In summary, the following cost function is used in this paper:
\begin{equation}
    \min \int_{t_0}^{t_f} \frac{\exp\left(\exp(N\omega(t))\right)}{M}\ dt
\end{equation}
\section{Simulation Result}
\begin{figure}[h]
    \centering
    \includegraphics[width=0.42\textwidth]{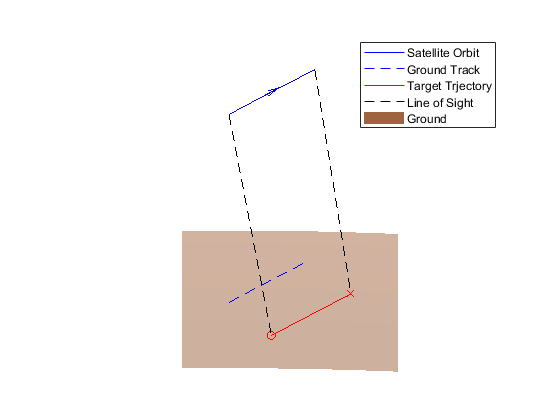}
    \includegraphics[width=0.42\textwidth]{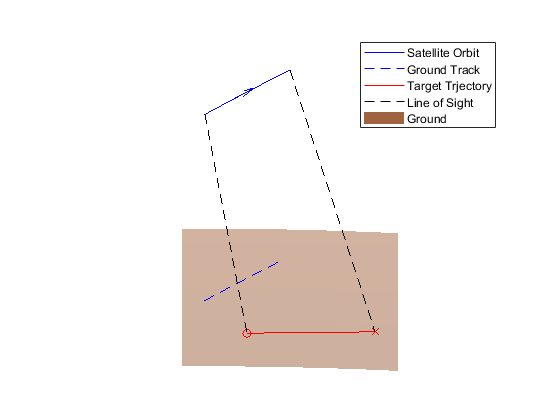}
    \includegraphics[width=0.42\textwidth]{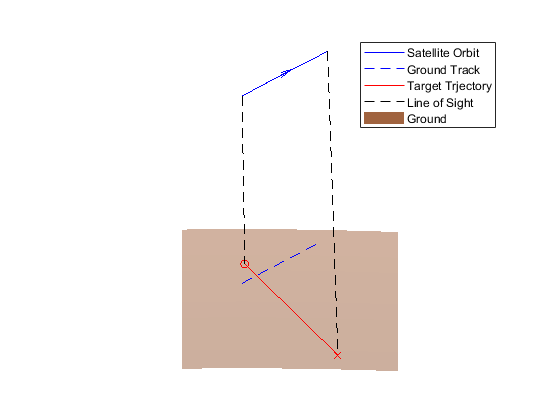}
    \includegraphics[width=0.42\textwidth]{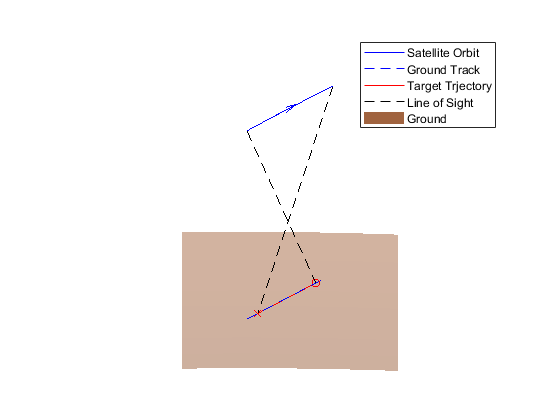}
    \caption{Four STRIP scenarios: (a) parallel, (b) offset, (c) perpendicular, (d) reverse}
    \label{fig:scenarios}
\end{figure}

The altitude of the satellite is set to 500 km, and the imaging duration is 30 seconds for each scenario. We have considered four different scenarios: parallel STRIP, offset STRIP, perpendicular STRIP, and reverse STRIP. The orbital trajectory and ground profile for each scenario are illustrated in Figure.~\ref{fig:scenarios}. The blue solid line represents the orbit, the red solid line shows the ground profile, and the dashed black lines are line of sight at the beginning and end of the imaging operation. The red circle marker and cross marker indicate the starting and end position of the STRIP, respectively.

For every scenario, the following initial guesses are used with time step of $1$ second:
\begin{equation} \label{eq:linear_sol}
    \bar s(t) = s(t_0) + \frac{s(t_f) - s(t_0)}{t_f - t_0} t, \quad
    \bar u(t) = \frac{s(t_f) - s(t_0)}{t_f - t_0}, \quad
    \bar \nu = 0
\end{equation}
the following initial Lagrangian and penalty weights are used:
\begin{equation}
    \mu = 1,\quad \lambda = 0
\end{equation}
and the following parameters are used:
\begin{equation}
    k_\nu = 0.5,\quad k_u = 0.5,\quad \gamma = 1.1, \quad \epsilon_V = \epsilon_g = 10^{-6}
\end{equation}
are used. The optimized profiles are compared with the linear solution(Linear) expressed in Eq.~\eqref{eq:linear_sol}.

\subsection{Optimal Scan Profile without $f_\text{CCD}$ constraints}
\begin{figure}[H]
    \centering
    \includegraphics[width=0.4\textwidth]{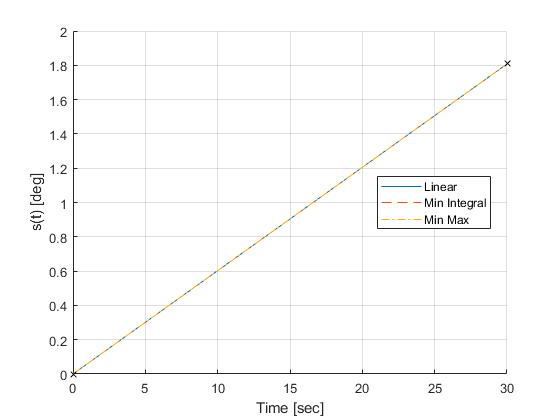}
    \includegraphics[width=0.4\textwidth]{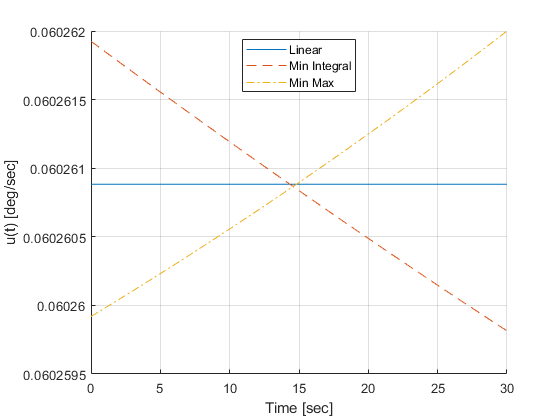}
    \includegraphics[width=0.4\textwidth]{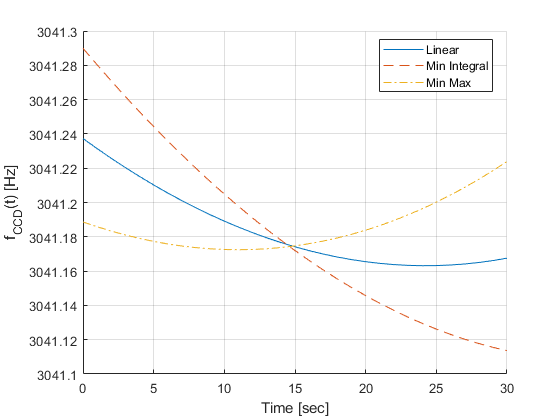}
    \includegraphics[width=0.4\textwidth]{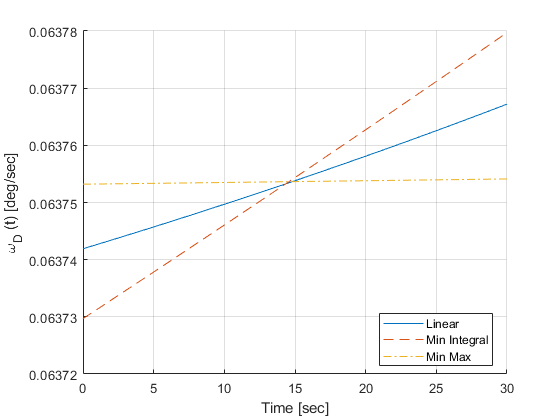}
    \caption{Simulation Result for Scenario 1: (a) $s(t)$, (b) $u(t)$, (c) $f_\text{CCD}$, (d) $\lVert \boldsymbol{\omega}_\mathcal{D} \rVert$}
    \label{fig:result_scn1}
\end{figure}

Firstly, we will show the optimized result without the $f_\text{CCD}$ constraints ($f_\text{LB} \leq f_\text{CCD}(t) \leq f_\text{UB}$). Table 1 summarizes the cost metric and feasibility of the terminal constraint for each scan profile in each scenario. Based on the result, all the optimized solutions satisfy the terminal constraint with a specified tolerance. Also, the minimizing integral solution(Min Integral) achieves minimum $\int \omega^2\ dt$, and minimizing maximum solution(Min Max) gets minimum $\text{max}\ \omega$ for all scenarios as expected. 
\begin{table}[H]
\caption{Results Summary of Scenarios \label{tab:case_table}}
\centering
\resizebox{\columnwidth}{!}{%
    \begin{tabular}{l cccc cccc ccc} \hline
        & \multicolumn{3}{c}{$\int \omega^2\ dt$} & & \multicolumn{3}{c}{$\max \omega(t)$}
        & & \multicolumn{3}{c}{$\lvert s(t_f) - s_f \rvert$} \\
        \cline{2-4} \cline{6-8}  \cline{10-12} 
    Scenario & Linear & Min Integral & Min Max && Linear & Min Integral & Min Max && Linear & Min Integral & Min Max \\
    \hline
    1  &  0.121936 & \textbf{0.121936} & 0.121936 && 0.063767 & 0.063780 & \textbf{0.063754} && \textbf{0} & 1.038e-09 & 3.444e-11\\
    2  &  8.313901 & \textbf{8.294143} & 8.339831 && 0.571981 & 0.553300 & \textbf{0.527831} && \textbf{0} & 5.555e-10 & 1.531e-10\\
    3  &  65.99975 & \textbf{65.01489} & 65.10274 && 1.637640 & 1.485586 & \textbf{1.473794} && \textbf{0} & 3.883e-09 & 1.544e-09\\
    4  &  58.62090 & \textbf{58.52228} & 58.52391 && 1.461166 & 1.410016 & \textbf{1.397096} && \textbf{0} & 9.673e-09 & 2.313e-09\\
    \hline
    \end{tabular}
}
\end{table}
The first set of figures in Figure.~\ref{fig:result_scn1} presents the results of scenario 1. Each figure shows the scan profile ($s(t)$), control command ($u(t)$), the scan frequency ($f_\text{CCD}$), and norm of the angular velocity command ($\lVert\boldsymbol{\omega}_\mathcal{D} \rVert $). Based on the result, the optimized results are not that different from the linear solution. One can see that the scan rate remains almost constant for all methods, and this demonstrates why TDI camera with fixed $f_\text{CCD}$ is limited to the parallel STRIP operation. Note that the target profile used in this scenario is almost parallel to the orbit but not exactly, resulting in a linearly increasing angular rate profile for the linear solution, although the deviation is practically negligible.  

\begin{figure}[H]
    \centering
    \includegraphics[width=0.4\textwidth]{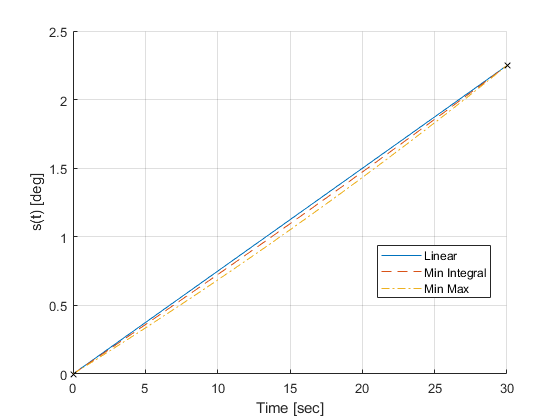}
    \includegraphics[width=0.4\textwidth]{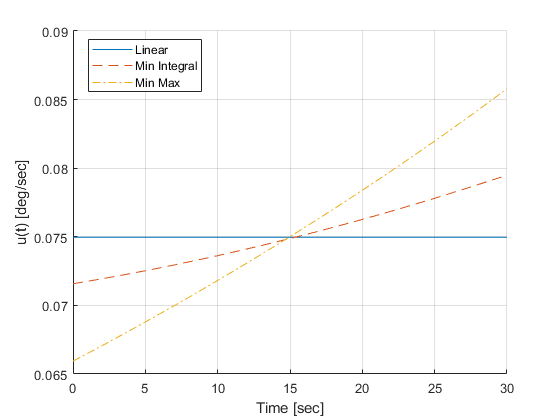}
    \includegraphics[width=0.4\textwidth]{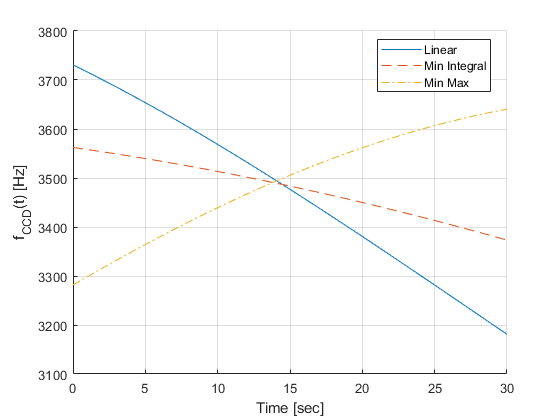}
    \includegraphics[width=0.4\textwidth]{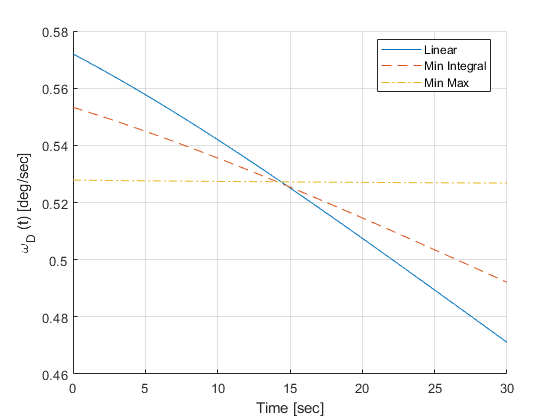}
    \caption{Simulation Result for Scenario 2: (a) $s(t)$, (b) $u(t)$, (c) $f_\text{CCD}$, (d) $\lVert \boldsymbol{\omega}_\mathcal{D} \rVert$}
    \label{fig:result_scn2}
\end{figure}

The second set of figures in Figure~\ref{fig:result_scn2} presents the results of scenario 2. Unlike the previous results, each method shows a distinct scan profile and angular velocity. Since the line of sight is small initially, as shown in Figure~\ref{fig:scenarios}, the required angular velocity is large at the beginning, and it decreases over time under the linear solution. By optimizing the solution, we can shape the scan profile. Notably, the angular velocity of the Min Max solution remains constant, while that of the Min Integral solution has an intermediate value.

\begin{figure}[h]
    \centering
    \includegraphics[width=0.4\textwidth]{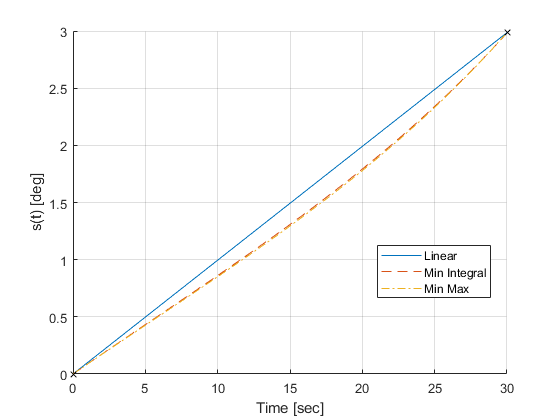}
    \includegraphics[width=0.4\textwidth]{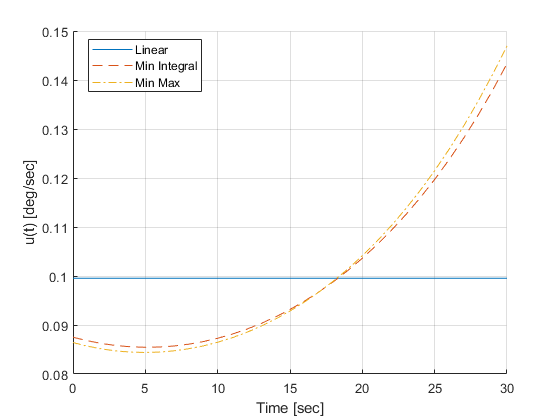}
    \includegraphics[width=0.4\textwidth]{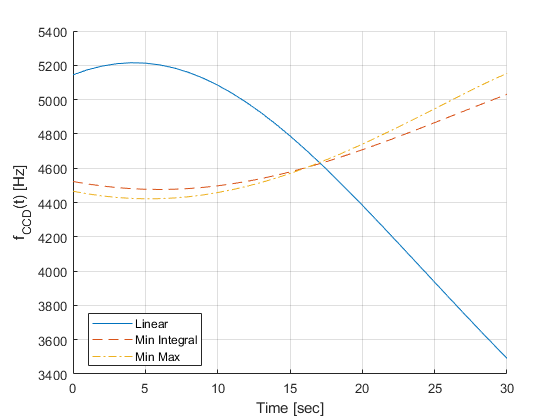}
    \includegraphics[width=0.4\textwidth]{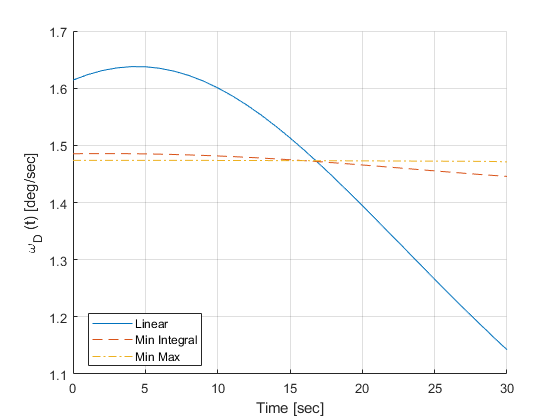}
    \caption{Simulation Result for Scenario 3: (a) $s(t)$, (b) $u(t)$, (c) $f_\text{CCD}$, (d) $\lVert \boldsymbol{\omega}_\mathcal{D} \rVert$}
    \label{fig:result_scn3}
\end{figure}

\begin{figure}[H]
    \centering
    \includegraphics[width=0.4\textwidth]{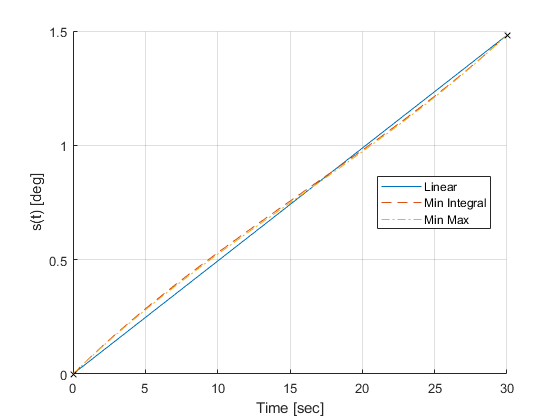}
    \includegraphics[width=0.4\textwidth]{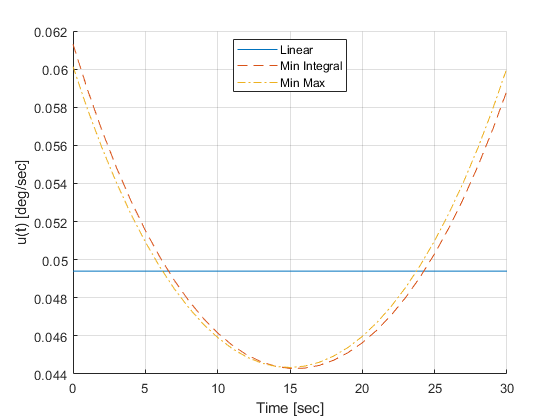}
    \includegraphics[width=0.4\textwidth]{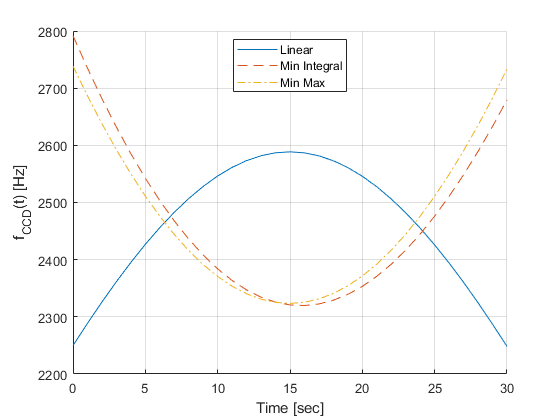}
    \includegraphics[width=0.4\textwidth]{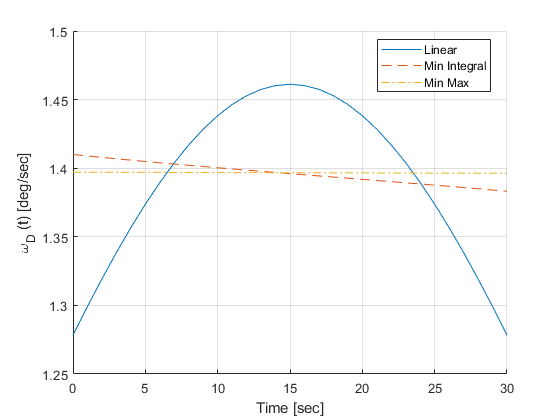}
    \caption{Simulation Result for Scenario 4: (a) $s(t)$, (b) $u(t)$, (c) $f_\text{CCD}$, (d) $\lVert \boldsymbol{\omega}_\mathcal{D} \rVert$}
    \label{fig:result_scn4}
\end{figure}
The results of scenarios 3 and 4 are shown in Figure.~\ref{fig:result_scn3} and~\ref{fig:result_scn4}, respectively. Both scenarios require high angular velocity due to the twisted relative motion between the satellite and the ground target. The angular velocity changes significantly under the linear solution, but the variations are reduced for both the Min Integral and Min Max solutions. 

Based on the simulation results, minimizing the integral of the norm is not very effective, as only a small portion is reduced after optimization. However, the Min Integral solution still has a preferred profile over the linear solution since it reduces the peak of angular velocity. On the other hand, minimizing the maximum norm can reduce the cost metric by up to 10 percent. The downside of the Min Max solution is that the DDP formulation becomes more sensitive due to the exponential cost term.

\subsection{Optimal Scan Profile with $f_\text{CCD}$ constraints}
In this subsection, the optimal scan profile for scenarios 3 and 4 are regenerated with the consideration of $f_\text{CCD}$ constraints ($f_\text{LB} \leq f_\text{CCD}(t) \leq f_\text{UB}$). The lower bound of $f_\text{CCD}$ is set as $4300$ and $2300$ Hz, and the upper bound is set as $4700$ and $2600$ Hz, respectively for each scenario. The profiles of the linear solution(blue), unconstrained minimum integral solution(red), and constrained minimum integral solution(yellow) are compared. 

\begin{figure}[h]
    \centering
    \includegraphics[width=0.43\textwidth]{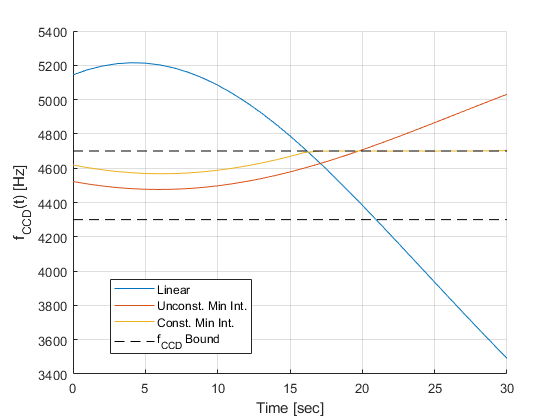}
    \includegraphics[width=0.43\textwidth]{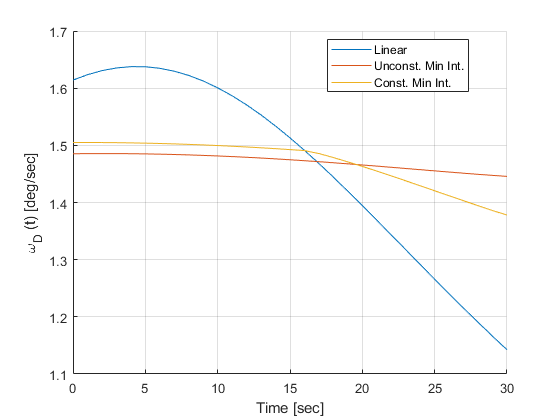}
    \caption{Simulation Result for Scenario 3 with $f_\text{CCD}$ constraint: (a) $f_\text{CCD}$, (b) $\lVert \boldsymbol{\omega}_\mathcal{D} \rVert$}
    \label{fig:result_const_scn3}
\end{figure}

\begin{figure}[h]
    \centering
    \includegraphics[width=0.43\textwidth]{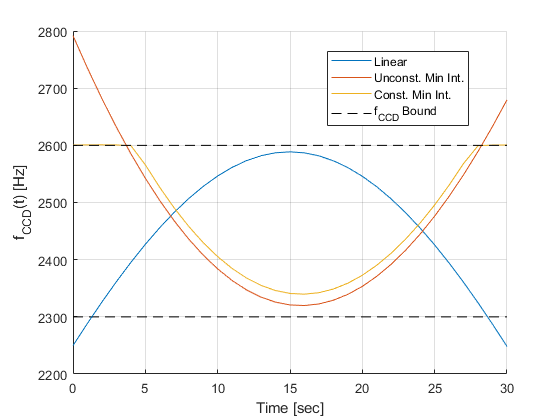}
    \includegraphics[width=0.43\textwidth]{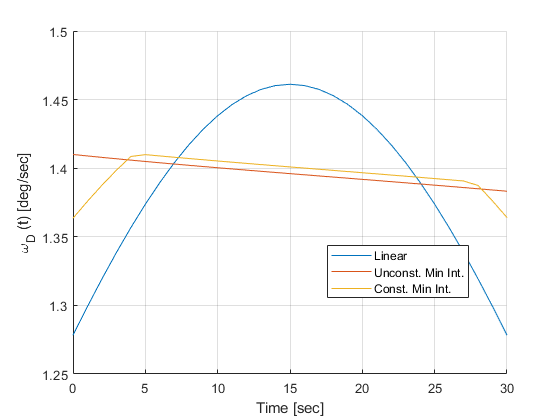}
    \caption{Simulation Result for Scenario 4 with $f_\text{CCD}$ constraint:  (a) $f_\text{CCD}$, (b) $\lVert \boldsymbol{\omega}_\mathcal{D} \rVert$}
    \label{fig:result_const_scn4}
\end{figure}

The results of scenarios 3 and 4 are shown in Figure.~\ref{fig:result_const_scn3} and~\ref{fig:result_const_scn4}, respectively. The results of the linear solution and Min Integral solution are identical to those of the previous subsection, and they violate the $f_\text{CCD}$ constraints. On the other hand, the constrained minimum integral solution satisfies the $f_\text{CCD}$ constraints which are represented as black dashed lines. That is, the path constraints are effectively managed by the augmented Lagrangian method. 

\begin{figure}[H]
    \centering
    \includegraphics[width=0.43\textwidth]{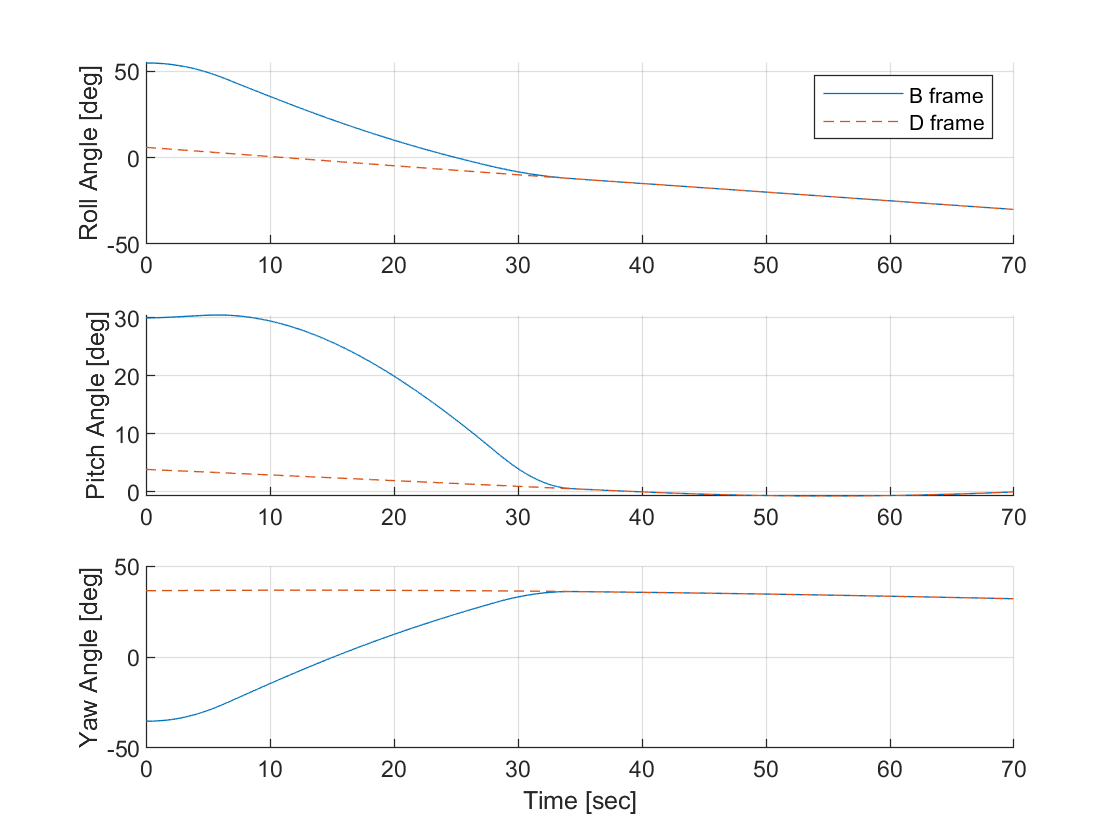}
    \includegraphics[width=0.43\textwidth]{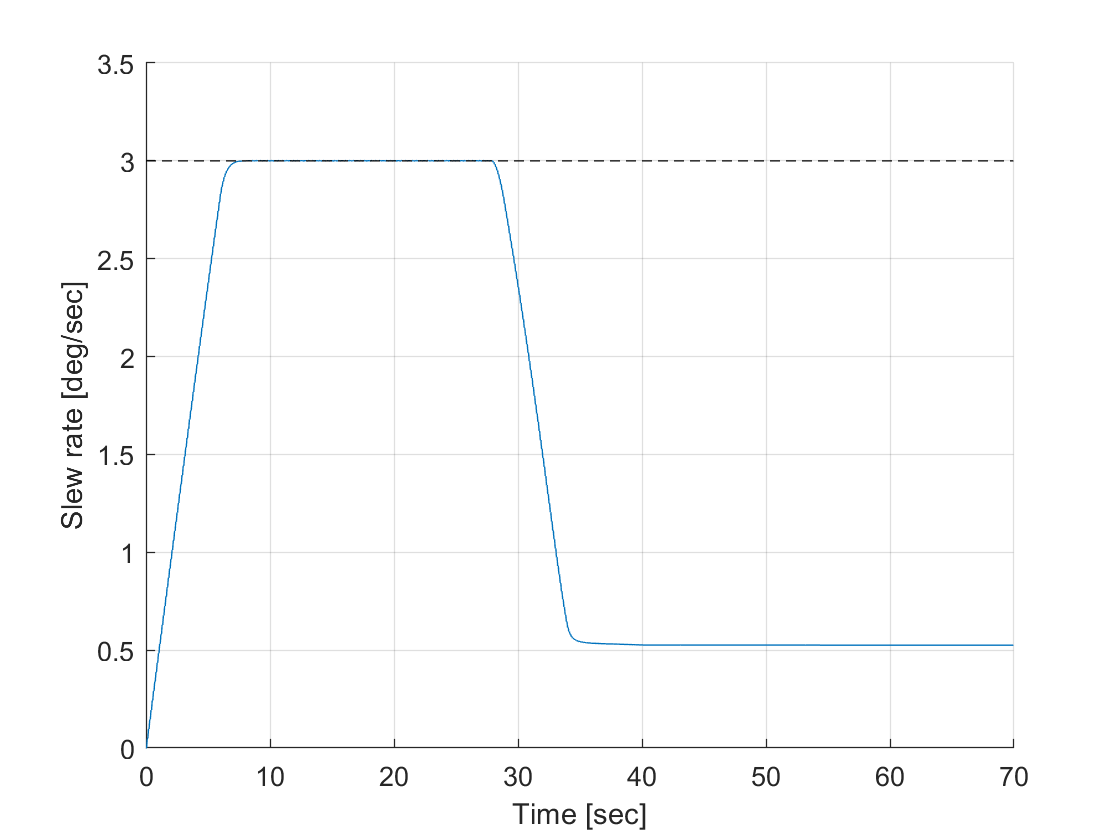}
    \caption{Simulation Result for Case 2:  (a) Euler 3-2-1 Angles of $\mathcal{B}$ and $\mathcal{D}$, (b) $\lVert \boldsymbol{\omega}_\mathcal{B} \rVert$}
    \label{fig:result_sim}
\end{figure}

\subsection{Implementation with Attitude Control Algorithm}
The Min Max $\omega_\mathcal{D}$ profile for scenario 2 is tested jointly with the attitude control algorithm. The attitude control algorithm and satellite parameters are adopted from the work to test the generated optimal profile\cite{han2024control}. The time step is adjusted to 0.1 seconds to match the 10 Hz control frequency, and the start time of the imaging operation is shifted by 40 seconds to allow for maneuvering time. 

The simulation results are shown in Figure.~\ref{fig:result_sim}, with the attitude profile represented in Euler angles and the norm of the body frame angular velocity displayed. Based on the simulation results, the satellite reorients its attitude from 0 to 40 seconds and maintains the strip imaging attitude thereafter. The slew rate of the satellite $\omega_\mathcal{B}$ becomes identical to the profile of $\omega_\mathcal{D}$ after 40 seconds (compare Min Max $\omega_\mathcal{D}$ of Figure.\ref{fig:result_scn2} and $\omega_\mathcal{B}$ of Figure.~\ref{fig:result_sim}). These results demonstrate that the generated profile is feasible and compatible with the general attitude-tracking control algorithm.

\section{Conclusion}
In this paper, we have derived the exact attitude command satisfying the TDI imaging constraints, which are scan rate and scan angle. The additional degree of freedom on the scan rate allows us to formulate an OCP, and we have considered two different cost metrics. The OCP is then solved via CDDP effectively, and the results are demonstrated through various practical scenarios.

Although the DDP works well for the simulation cases, the algorithm is highly sensitive to the initial guess and parameters, requiring several initial trials before achieving a good result. Once proper parameters are found, the method consistently provides a converged solution. On the other hand, the method works much faster than conventional optimization methods.

Lastly, minimizing the integral of the rate squared turns out to be ineffective, while minimizing the maximum rate is effective. However, the integral-minimizing solution still has a preferred profile over the linear solution, and its DDP formulation is less sensitive to the initial parameters compared to the maximum-minimizing solution.

\bibliographystyle{AAS_publication}   
\bibliography{references}   

\end{document}